\documentclass[graybox]{svmult}

\usepackage{mathptmx}       

\usepackage{helvet}         
\usepackage{courier}        

\usepackage{graphicx}              

\usepackage{siunitx}
\usepackage{amsbsy}
\usepackage{amsmath,bm}
\usepackage{amssymb}
\usepackage{amsfonts}
\usepackage{algorithmic}
\usepackage{algorithm} 

\usepackage{multirow}

\usepackage[bottom]{footmisc}
\usepackage{cite}

\usepackage{url}
\usepackage{xcolor}

\newcommand{\R}{\mathbb R}
\newcommand{\N}{\mathbb N}

\renewcommand{\rho}{\varrho}

\begin{document}

\title*{Operator splitting for port-Hamiltonian systems}
\titlerunning{Operator splitting for  port-Hamiltonian systems}
\author{Andreas Frommer
\and 
Michael Günther 
\and
Björn Liljegren-Sailer  
\and \\
Nicole Marheineke  
}
\institute{Andreas Frommer and Michael Günther \at Bergische Universität Wuppertal, IMACM, Gaußstraße 20, D-42119 Wuppertal,\\
\email{[frommer,guenther]@uni-wuppertal.de}
\\
Björn Liljegren-Sailer and Nicole Marheineke \at Universität Trier, FB IV -- Mathematik, Universitätsring 15,
D-54296 Trier, \\ \email{[bjoern.sailer,marheineke]@uni-trier.de}}

%
%
\maketitle

\abstract{The port-Hamiltonian approach presents an energy-based modeling of dynamical systems with energy-conservative and energy-dissipative parts as well as an interconnection over the so-called ports. In this paper, we apply an operator splitting that treats the energy-conservative and energy-dissipative parts separately. This paves the way for linear equation solvers to exploit the respective special structures of the iteration matrices as well as the  multirate potential in the different right-hand sides. We illustrate the approach using test examples from coupled multibody system dynamics.}

  \section{Introduction}

Operator splitting~\cite{mclachlan_quispel_2002} is an efficient tool to numerically solve initial-value problems of ordinary differential equations (ODE-IVPs)
  \begin{equation}
  \label{ode}
      \dot x= f(x,t)=f_1(x,t)+f_2(x,t), \qquad x(t_0)=x_0,
  \end{equation}
that allow for a splitting of the right hand side $f(x,t)$ into two parts $f_1(x,t)$ and $f_2(x,t)$ of profoundly different behaviour, e.g., with respect to stiffness, computational costs, dynamics etc.
Rewriting the system in the homogeneous from
  \begin{eqnarray}
  \begin{pmatrix}
      \dot x \\ \dot s
  \end{pmatrix}
  & = & 
  \begin{pmatrix}
  f(x,s) \\ 1
  \end{pmatrix} =
  \underbrace{\begin{pmatrix}
      f_1(x,s) \\
      1
      \end{pmatrix}}_{\displaystyle \tilde f_1(x,s):=}
      +
      \underbrace{\begin{pmatrix}
      f_2(x,s) \\
      0
      \end{pmatrix}}_{\displaystyle \tilde f_2(x,s):=}
      , \quad
      \begin{pmatrix}
          x(t_0) \\ s(t_0)
      \end{pmatrix}
      =
      \begin{pmatrix}
          x_0 \\ 0
      \end{pmatrix},
  \end{eqnarray}
  the idea is to alternately solve dynamical systems driven by $\tilde f_1$ and $\tilde f_2$, respectively.
    
Strang splitting~\cite{strang1968} 
    solves the ODE system with respect to the first, second and again first split right-hand side with step sizes $h/2$, $h$ and {$h/2$}, where the initial values are given by the respective final values, i.e.,
  \begin{align}
      \begin{pmatrix}
          \dot x_{1/3} \\ \dot s_{1/3}
      \end{pmatrix}& = \tilde f_1(x_{1/3},s_{1/3}),
      \quad && \begin{pmatrix}
          x_{1/3}(t_0) \\ s_{1/3} (t_0) \end{pmatrix} =
          \begin{pmatrix} x_0 \\ 0
      \end{pmatrix}, \label{strang.1} \\
\begin{pmatrix}
          \dot x_{2/3} \\ \dot s_{2/3}
      \end{pmatrix}& = \tilde f_2(x_{2/3},s_{2/3})
      \quad &&\begin{pmatrix}
          x_{2/3}(t_0) \\ s_{2/3} (t_0) \end{pmatrix} =
          \begin{pmatrix} x_{1/3}(h/2) \\ h/2
      \end{pmatrix}, \label{strang.2} \\
        \begin{pmatrix}
          \dot x_{1} \\ \dot s_{1}
      \end{pmatrix}& = 
      \tilde f_{1}(x_1,s_1) \quad &&\begin{pmatrix}
          x_{1}(t_0) \\ s_{1}(t_0) \end{pmatrix} =
          \begin{pmatrix} x_{2/3}(h) \\ h
      \end{pmatrix}. \label{strang.3} 
  \end{align}
  This scheme provides an order two approximation $x_1(h/2)$ for the exact solution $\varphi(h) \circ \begin{pmatrix} 
  x_0 \\ 0 \end{pmatrix}$ of~\eqref{ode} at time point $t_0+h$, starting at the initial value $x_0$ at $t_0$. In short hand:
  \begin{equation}
      \begin{pmatrix}
          x_1(h/2) \\ s_1(h/2)
      \end{pmatrix} =  \varphi_{\tilde f_1}(\frac{h}{2}) \circ
      \varphi_{\tilde f_2}(h) \circ \varphi_{\tilde f_1}(\frac{h}{2}) \circ \begin{pmatrix}
          x_0 \\ 0
      \end{pmatrix} = 
      \varphi(h) \circ \begin{pmatrix}
          x_0 \\ 0
      \end{pmatrix}+  \begin{pmatrix}
      {\mathcal{O}}(h^3)\\ 0
      \end{pmatrix}
  \end{equation}
  where $\varphi_{\tilde f_1}(h)$ and  $\varphi_{\tilde f_2}(h)$ denote the solution of the flow $\tilde f_1$ and $\tilde f_2$, respectively, at time point $t_0+h$, starting at $t_0$. 

If the three initial-value problems~\eqref{strang.1}-\eqref{strang.3} are solved by consistent approximation schemes $\Phi_{\tilde f_1}(h/2)$, $\Phi_{\tilde f_2}(h)$ and $\Phi_{\tilde f_1}(h/2)$ with step sizes $h/2$, $h$ and {$h/2$}, the corresponding numerical approximation 
$$\tilde \Phi(h) \circ \begin{pmatrix} x_0 \\ 0 \end{pmatrix}:= \Phi_{\tilde f_1}(h/2) \circ\Phi_{\tilde f_2}(h) \circ\Phi_{\tilde f_1} (h/2) \circ \begin{pmatrix} x_0 \\ 0 \end{pmatrix}$$ is of second order due to the symmetry of the approach. We will apply this idea of operator splitting to the special case of port-Hamiltonian ODE systems.

  The paper is organized as follows: Port-Hamiltonian ODE systems are brief\/ly introduced in Section 2. The following sections shows how the operator splitting approach can be applied to port-Hamiltonian systems, with discrete gradient schemes being the method-of-choice.
  Sections 5 and 6 show how this approach can exploit both the 
  multirate behaviour and 
  system structure at the level of nested integration and linear solvers. Numerical results for test example from coupled multibody systems are discussed in Section 7. Finally, we give a summary and outlock to future work.

  \section{Port-Hamiltonian ODE systems}
  
  An ODE-IVP of port-Hamiltonian structure~\cite{vanderschaft2006} is given by
  \begin{align}
  \label{port-ham.a}
      \dot x & = (J(x)-R(x)) \nabla H(x) + B(x) u(t), \\
      \label{port-ham.b}
      y & = B(x)^\top \nabla H(x)
  \end{align}
  with state variable $x: [t_0,t_{end}] \rightarrow \R^n$, a twice continuously differentiable Hamiltonian $H: \R^n \rightarrow \R$, $n \times n$ matrices $J(x)=-J(x)^\top$ skew-symmetric and $R(x)\ge 0$ positive semi-definite, an input signal $u: [t_0,t_{end}] \rightarrow \R^d$, an output $y: [t_0,t_{end}] \rightarrow \R^d$ and a $n \times d$-dimensional port matrix $B(x)$.
  
  A fundamental property for the solution $x(t)$ of a port-Hamiltonian system is the dissipativity inequality
  \begin{align*}
  \frac{d}{dt} H(x(t)) &=  \nabla H(x(t))^\top \dot x(t)  
    = - \nabla H(x(t))^\top R(x(t)) \nabla H(x(t)) + y(t)^\top u(t) \\ \nonumber &\le y(t)^\top u(t).
    \end{align*}
This can also be written in the integral form
    \begin{align} \label{diss.ineq}
      H(x(t_0+h))&-H(x(t_0)) = \nonumber \\
      & = - \int_{t_0}^{t_0+h} \nabla H(x(s))^\top R(x(s)) \nabla H(x(s)) \, ds   
       + \int_{t_0}^{t_0+h} y(s)^\top u(s) \, ds \\      
       & \le   \int_{t_0}^{t_0+h} y(s)^\top u(s) \, ds\nonumber 
  \end{align}
  Without input and dissipation $R \equiv 0$, the Hamiltonian is an invariant of the port-Hamiltonian ODE~\eqref{port-ham.a}.

If port-Hamiltonian systems are solved numerically, the aim is not only efficiency, but also structure preservation. We aim for schemes that preserve the dissipativity inequality~\eqref{diss.ineq} at a discrete level. 

\section{Operator splitting for port-Hamiltonian systems}
The splitting of the right-hand side of the port-Hamiltonian ODE~\eqref{port-ham.a} into an energy-preserving part $f_1(x,t):=J(x) \nabla H(x)$ and a dissipative energy-coupling part $f_2(x,t):=-R(x) \nabla H(x) + B(x) u(t)$ is a natural choice.

The approximation $z(h)$, given by the Strang splitting approximation 
\begin{align}
\label{phs.strang.1}
\begin{pmatrix} \dot v \\ \dot s_v \end{pmatrix} & =  \begin{pmatrix} - R(v) \nabla H(v)+ B(v) u(s_v) \\ 1 \end{pmatrix}, \qquad &&\begin{pmatrix} v(t_0) \\ s_v(t_0) \end{pmatrix} = \begin{pmatrix} x_0 \\ 0 \end{pmatrix}, \\
\label{phs.strang.2}
\begin{pmatrix} \dot w \\ \dot s_w \end{pmatrix}  & =  \begin{pmatrix} J(w) \nabla H(w) \\ 0 \end{pmatrix}, \qquad &&\begin{pmatrix} w(0) \\ s_w(t_0) \end{pmatrix} = \begin{pmatrix} v(h/2) \\ h/2 \end{pmatrix}, \\
\label{phs.strang.3}
\begin{pmatrix} \dot z \\ \dot s_z \end{pmatrix}  & =  \begin{pmatrix} - R(z) \nabla H(z) + B(z) u(s_z) \\ 1 \end{pmatrix}, \qquad && \begin{pmatrix} z(t_0) \\ s_z(t_0) \end{pmatrix} = \begin{pmatrix} w(h) \\ h \end{pmatrix} \\
y_z & =   B(z)^\top \nabla H(z), \nonumber 
\end{align}
 provides an order two approximation $z(h/2)$
 for the exact solution $\varphi(h) \circ \begin{pmatrix} x_0 \\ 0 \end{pmatrix}$
 of~\eqref{port-ham.a} 
  at time point $t_0+h$ and the output $y$, respectively, starting at the initial value $x_0$ at $t_0$. In short hand:
 \begin{eqnarray*}
      \begin{pmatrix}
          z(h/2) \\ s_z(h/2)
      \end{pmatrix} &= & \varphi_{\tilde f_1}(\frac{h}{2}) \circ
      \varphi_{\tilde f_2}(h) \circ \varphi_{\tilde f_1}(\frac{h}{2}) \circ \begin{pmatrix}
          x_0 \\ 0
      \end{pmatrix} = 
      \varphi(h) \circ \begin{pmatrix}
          x_0 \\ 0
      \end{pmatrix}+  \begin{pmatrix}
      {\mathcal{O}}(h^3)\\ 0
      \end{pmatrix} \\
      y_z(h/2) & = & y(h) +{\mathcal{O}}(h^3)
  \end{eqnarray*}
 %
%
The dissipativity inequality is fulfilled exactly at a discrete level, i.e.,
\begin{align*}
    H(z(h/2))-H(x_0) & = (H(z(h/2))-H(z(0)))+(H(w(h))-H(w(0)))+ \\
    & \quad + (H(v(h/2))-H(x_0)) \\
    & \le 
      \int_0^{h/2} y_v(\tau)^\top u(\tau) \, d\tau +
      \int_0^{h/2} y_z(\tau)^\top u(\tau+h/2) \, d\tau.
\end{align*}

%

\section{Discrete gradient methods}
Applying operator splitting to the port-Hamiltonian system~\eqref {port-ham.a}, \eqref{port-ham.b} requires the following demands:
\begin{itemize}
    \item the numerical approximations of~\eqref{phs.strang.1} and \eqref{phs.strang.3} 
    are dissipative at a discrete level for vanishing $u \equiv 0$, i.e., it holds in this case
     \begin{align*}
      H(v(h/2))  \leq  H(x_0) \hspace{0.8cm} \text{and} \hspace{0.8cm} {H(z(h/2))  \leq  H(w(h)) } 
      \end{align*}
      holds for all step sizes $h>0$. 
    \item the numerical approximation $w(h)$ 
    of~\eqref{phs.strang.2} is energy preserving at a discrete level, i.e., 
    \begin{align*}
     {H(w(h))  =  H(v(h/2)) }  
      \end{align*}
      holds for all step sizes $h>0$.
\end{itemize}
Though~\eqref{phs.strang.2} defines a symplectic flow for $J$ regular, symplectic schemes are only energy preserving for quadratic Hamiltonians, but not in general.
Hence methods of choice for solving~\eqref{port-ham.a}-\eqref{port-ham.b} here are discrete gradient methods~\cite{GM_Go96}, defined by
\begin{align}
    \frac{x_1-x_0}{h} & = (\bar J(x_0,x_1,h)- \bar R(x_0,x_1,h)) \bar \nabla H(x_0,x_1) + \bar B(x_0,x_1,h) \bar u(t_0,t_1,h), \label{dgm.1}  \\
    \label{dgm.2}
    y_1 & = \bar B(x_0,x_1,h)^\top \bar \nabla H(x_0,x_1),
\end{align}
with the discrete gradient $ \bar \nabla H: \R^n \times \R^n \rightarrow \R$ fulfilling the two conditions
\begin{align*}
   i)& \hspace{0.5cm} \bar \nabla H(x,y)^\top (y-x)  = H(y)-H(x), \\
    ii)& \hspace{0.5cm} \bar \nabla H(x,x)  = \nabla H(x)
\end{align*}
for all $x,y \in \R^n$, and the skew-symmetric matrices $\bar J(x,y,h)$, the positive-semidefinite matrices $\bar R(x,y,h)$, the matrices $\bar B(x,y,h)$ and $\bar u(t_0,t_1,h)$ fulfilling the compatability conditions
$\bar J(x,x,0)=J(x), \bar R(x,x,0)=R(x), \bar B(x,x,h)=B(x)$ and $\bar u(t_0,t_0,0)=u(t_0)$.

By construction, the discrete gradient method satisfies a discrete version of the 
dissipativity inequality \cite{mclachlan1999geometric}:
\begin{eqnarray*}
H(x_1)-H(x_0) & = & -h \bar \nabla H(x_0,x_1,h)^\top \bar R(x_0,x_1,h) 
 \bar \nabla H(x_0,x_1,h)+ h y_1^\top \bar u(t_0,t_1,h).
 \end{eqnarray*}

One notes that discrete gradient methods are implicit and {limited to order two} for non-vanishing $R$,
whereas order one would have been sufficient for deriving an overall order two operator splitting scheme.

\section{Multirate potential} \label{sec:multirate}

A typical situation is that the energy preserving part defines a highly oscillatory behaviour, which requires small stepsizes, whereas the dissipative part may be characterized by a very slow dynamics, allowing for large step sizes.

One idea to exploit this multrate behaviour is to use low (order two) schemes for the dissipative part, as these are sufficient for obtaining accurate results also for large stepsizes, but to use highly accurate numerical schemes for the energy preserving part, which then also allows for larger step sizes.

Another way of exploiting this multirate behaviour given by different parts of the right-hand side with a time constants ratio of $m \in \N$ is to combine operator splitting with nested integration~\cite{SEXTON1992665}, i.e., to replace the step~\eqref{strang.2} with step size $h$ by $m$ subsequent steps of step size $h/m$:\\

\begin{quote}
\normalsize
\noindent
$\dot w_1  =  J(w_1) \nabla H(w_1), \quad w_1(0)=v(h/2)$ \\
for $i=2,\ldots,m:$ $\dot w_i  =  J(w) \nabla H(w), \qquad w_i(0)=w_{i-1}(h)$\\
$w(h):=w_m(h)$
\end{quote}

A straightforward way to derive higher order methods are composition methods. Using the operator splitting approach based on the discrete gradient method as the base scheme of order two, a composition of three base schemes with three different step sizes yields an order four scheme. However, composition schemes with order higher than two demand negative time steps~\cite{suzuki1991}, which contradicts the discrete dissipativity condition in the case of $R \neq 0$. 

In lattice quantom chromodynamics, where the gauge action plays the role of the fast and cheap part, and the fermionic action plays the role of the slow and expensive part this problem could be circumvented by the following idea:  if one uses an integrator of second
order for the slow action with step size $h$, and approximates the fast action by $m$ steps of the second order scheme with step size 
$h/(2m)$, the error of the overall multirate scheme will be of order ${\mathcal O} (h^2)+{\mathcal O}(( \frac{h}{m})^2)+
{\mathcal O}  (h^4)$. With the use of force gradient~\cite{OMELYAN2003272} information only at the slowest level it is possible
to cancel the leading error term of order ${\mathcal O (h^2)}$. If the mutirate factor $m$ between the time constants of both schemes is high enough, one gets already for quite small step size $h \leq \frac{1}{m}$  the overall order is then given by the leading error term of order ${\mathcal O} (h^4)$, i.e., the
scheme has an effective order of four~\cite{GM_sh18}. Whether this approach can be applied to the port-Hamiltonian setting is an open question.

Another idea is to use highly accurate numerical schemes for the energy preserving part, which allows for larger step sizes.

\section{Linear solvers}

Operator splitting also allows for exploiting the special structure one obtains when solving the energy-preserving and dissipative subsystems.
For simplicity of discussion, we assume constant matrices $J,R$ and $B$ and a quadratic Hamiltonian $H(x)=x^\top Q x$ here. Note that the discrete gradient method is equivalent to the implicit midpoint rule in case of linear systems.

{Applying the standard average vector field method as a classical discrete gradient method of order two, one yields the following linear system 
\begin{eqnarray} \label{eq:J_average_field}
    \frac{x_1-x_0}{h}  =  J Q \frac{x_1+x_0}{2} & \Rightarrow &
    \left(I -\frac{h}{2} \tilde J\right) \tilde x_1  =  \left(I + \frac{h}{2} \tilde J \right) \tilde x_0
\end{eqnarray}
for the energy preserving part~\eqref{phs.strang.2},
where we have applied a congruency with $Q^{1/2}$ with
$\tilde x_1 = Q^{1/2} x_1$, $\tilde x_0= Q^{1/2} x_0$ and
$\tilde J= Q^{1/2} J Q^{1/2}$.
}

{
Computing approximations  for $\tilde{x}_1$ in \eqref{eq:J_average_field} by applying a linear solver to the matrix $I -\frac{h}{2} \tilde J$ is not a structure-aware approach, since the computed approximations will, typically, not reflect energy preservation, at least if, for efficiency reasons, we do not aim at a very accurate solution. A structure aware approach to compute approximations for \eqref{eq:J_average_field} arises, if we use a {\em matrix function} approach, i.e.\ we consider
\[
\tilde{x}_1 = C(\tilde{J},\tfrac{h}{2})\tilde{x_0},
\]
where $C(\tilde{J},\frac{h}{2})$ is the matrix function evaluation of the Cayley transform
\[
C(z,a) = \frac{z+a}{z-a}, \enspace z \in \mathbb{C}\setminus \{a\}.
\]
We can now use the Arnoldi method \cite{FrSi06} to evaluate the action of $C(\tilde{J},\frac{h}{2})$ on the vector $\tilde{x}_0$. Each iteration of this method requires one matrix-vector multiplication with $\tilde{J}$, and unlike the general Arnoldi method it relies on a short recurrence because $\tilde{J}$ is skew-hermitian. This approach is structure preserving, since all iterates it produces will have the same 2-norm than $\tilde{x}_0$. 
This is because the Arnoldi method obtains its $k$-th iterate as the Cayley transform with a matrix which is the orthogonal projection of $\tilde{J}$ onto the $k$-th Krylov subspace. 
}

We can work in a similar manner with the original matrix $J$ by using the inner product defined by $Q$ rather than the Euclidian inner product; see \cite{ConGol76,FroKah2022,Gueetal2022}.

For the dissipative energy-coupling part~\eqref{phs.strang.1} one gets the linear system
\begin{eqnarray*}
    \frac{x_1-x_0}{h/2} & = & -R Q \frac{x_1+x_0}{2} + B \frac{u_(t_0+h/2)+u(t_0)}{2} \Rightarrow \\
    \left(I + \frac{h}{4} \tilde R \right) \tilde x_1 & = & \left(I -\frac{h}{4} \tilde R \right) \tilde x_0 + \tilde B \frac{u_(t_0+h/2)+u(t_0)}{2},
\end{eqnarray*}
where we have applied again a congruency with $Q^{1/2}$ with
$\tilde x_1 = Q^{1/2} x_1$, $\tilde x_0= Q^{1/2} x_0$ and
$\tilde R= Q^{1/2} R Q^{1/2}$, $\tilde B=Q^{1/2} B$. Now the iteration matrix $I + \frac{h}{2} \tilde R$ is symmetric positive-definite. {In a coupling context, $R$ will be composed of the dissipation operators of the individual systems, and we have the potential of using targeted preconditioners for each of the systems to obtain particulalry efficient solvers for the composed system, also for larger time steps. We do not go into further details here.}

\section{Numerical examples}
Finally we discuss two numerical examples: a two mass oscillator with damping~\cite{GM_gu22b} for testing the different time integration approaches, and a single mass-spring-damper (MSD) chain~\cite{GPBS2012} available from the Port Hamiltonian Benchmark System\footnote{\url{https://algopaul.github.io/PortHamiltonianBenchmarkSystems.jl}} to validate the structure preserving properties of the the matrix Arnoldi approach.

\subsection{A two mass osciallator with damping}
\begin{figure}[htb]
\begin{center}
\includegraphics[width=0.75\textwidth]{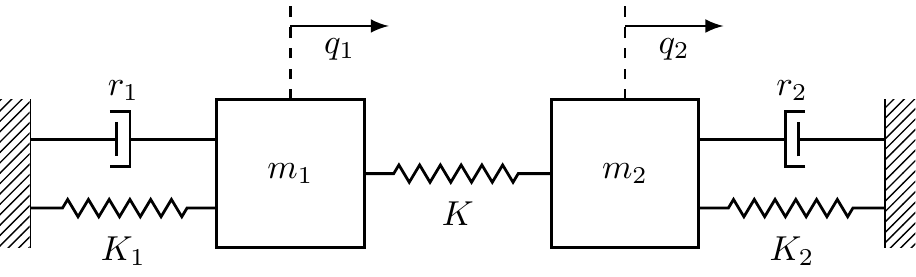}
\end{center}
\caption{\label{fig:two-masses_two-dampers_three-springs} ODE two masses oscillator with damping. The coordinates $q_1,\, q_2$ describe the position of the masses, taken from~\cite{GM_gu22b}.}
\end{figure}

As an example to validate the different numerical time integration approaches we consider two systems ($i=1,2$), each consisting of a mass $m_i>0$, which is connected  via a massless spring $K_i$ to walls with damping $r_i>0$, which are coupled by a spring $K$, see  Fig.~\ref{fig:two-masses_two-dampers_three-springs}.

\begin{figure}[htb]
\begin{center}
\includegraphics[width=\textwidth]{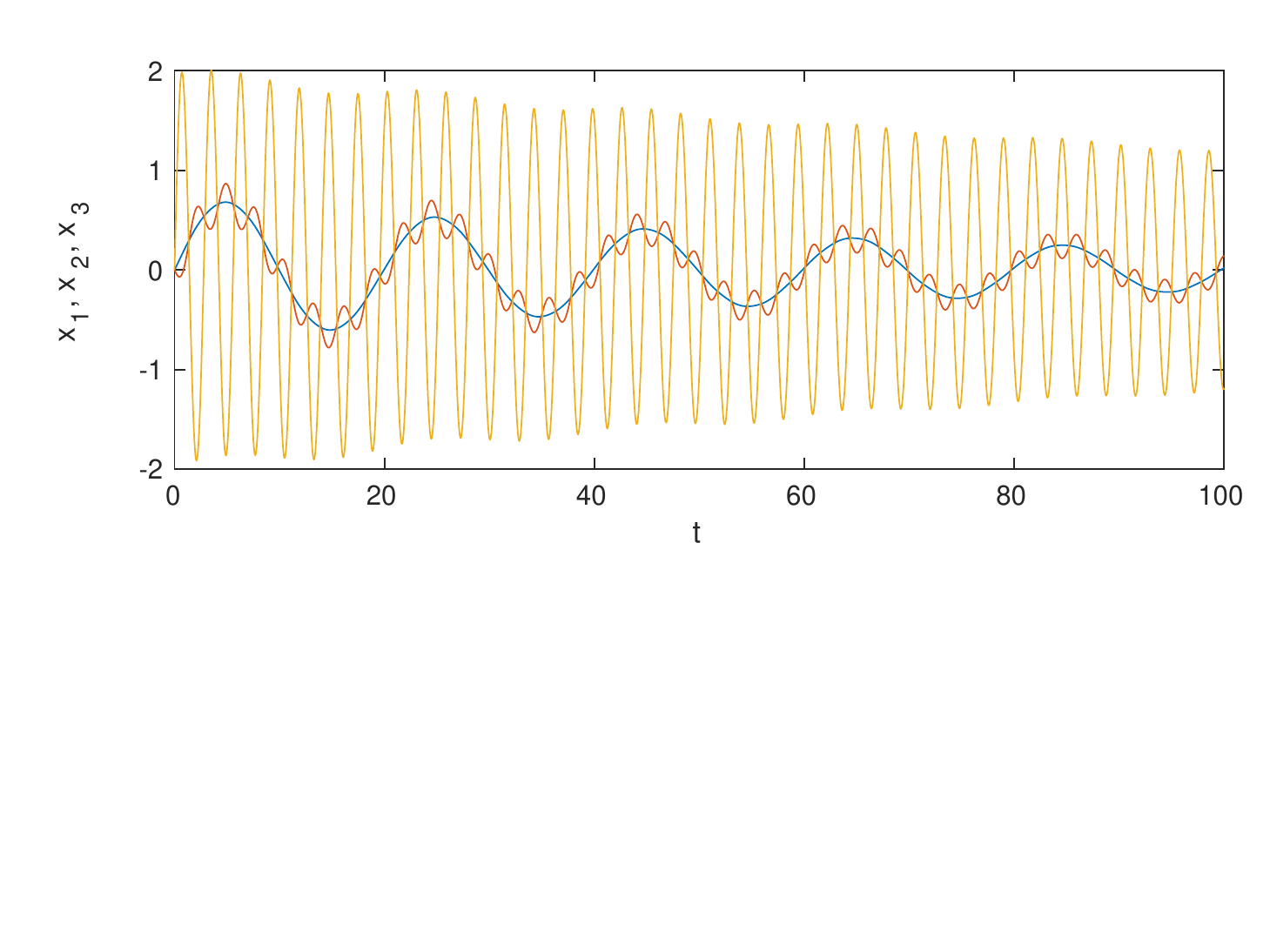}
\end{center}
\vspace*{-4cm}
\caption{\label{fig:solutions} Exact solution of the first three components of $x$ for the parameters $m_1=m_2=200, K=K_1=10, K_2=1000, r_1=5, r_2=2$}
\end{figure}

To set up the coupled system, we have for the two systems the positions $q_1,\, q_2$ and momenta $p_1,\, p_2$ for the masses $m_1,\, m_2$, respectively, as well as position $q=q_2$ as coupling variable between the two systems. The port-Hamiltonian description of this coupled system (with the coupling described by the off-diagonal elements of the skew-symmetric matrix $J$) is given by 
\begin{align*}
    \dot x & = (J-R) Q x, \qquad x(0)=x_0
\end{align*}
for the unknown 
$x:=(q_1,q_1-q,q_2,p_1,p_2)^\top$, with 
%
\begin{align*}
    & J= \begin{pmatrix}
    0 & 0 & 0 & 1 & 0 \\
    0 & 0 & 0 & 1 & -1 \\
    0 & 0 & 0 & 0 & 1 \\
   -1 & -1& 0 & 0 & 0 \\
    0 &1 &-1 & 0 & 0
    \end{pmatrix}, \quad
    R = \mbox{diag} (0,0,0,r_1,r_2), \quad
    Q = \mbox{diag} (K_1.K,K_2,\frac{1}{m_1},\frac{1}{m_2}),
\end{align*}
and the Hamiltonian given by $H(x)=x^\top Q x$. 
With $\tilde x:= Q^{1/2}x$ one gets 
\begin{equation}
    \dot {\tilde x}  =  (\tilde J- \tilde R) \tilde x
\end{equation}
with
\begin{align*}
    \tilde J =\begin{pmatrix}
    0 & 0 & 0 & \sqrt{\frac{K_1}{m_1}} & 0 \\
    0 & 0 & 0 & \sqrt{\frac{K}{m_1}} & -\sqrt{\frac{K}{m_2}}\\
    0 & 0 & 0 & 0 &  \sqrt{\frac{K_2}{m_2}}\\
   -\sqrt{\frac{K_1}{m_1}} & -\sqrt{\frac{K}{m_1}}& 0 & 0 & 0 \\
    0 &  \sqrt{\frac{K}{m_2}}&- \sqrt{\frac{K_2}{m_2}}& 0 & 0
    \end{pmatrix}, \quad
    R = \mbox{diag} \left(0,0,0,\frac{r_1}{m_1},\frac{r_2}{m_2}\right), 
\end{align*}
The waveforms of $x_1,x_2$ and $x_3$ are depicted in Fig.~\ref{fig:solutions}.

\begin{figure}[htb]
\begin{center}
\includegraphics[width=\textwidth]{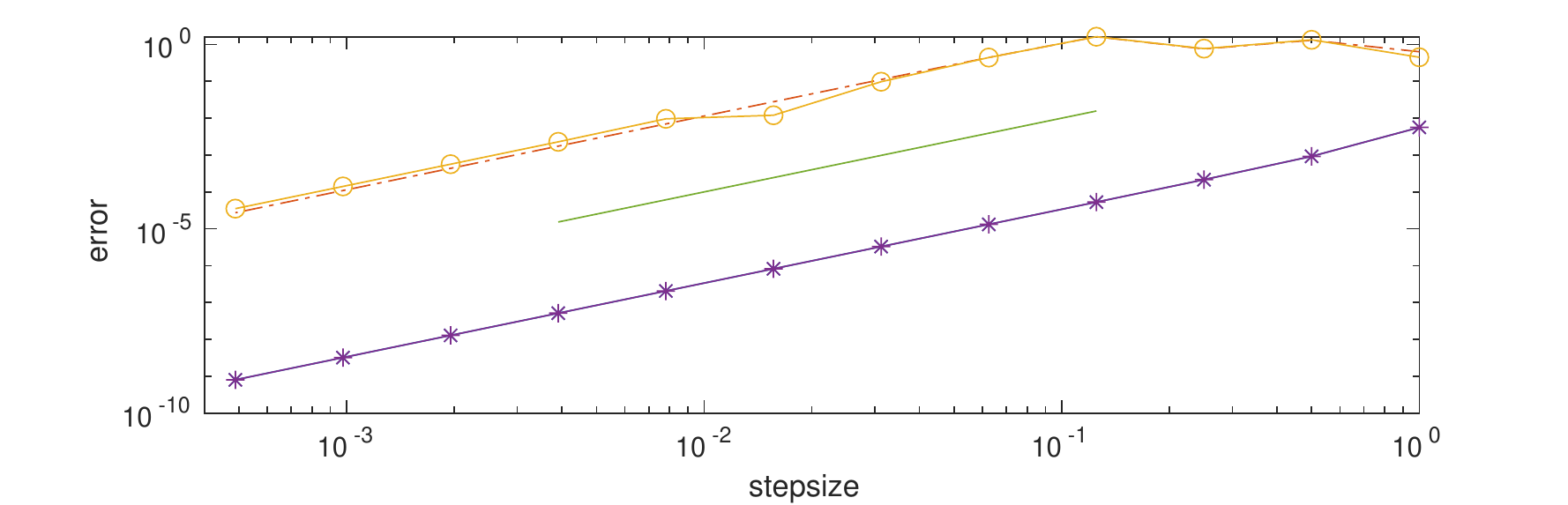}
\end{center}
\caption{\label{fig:solutions-2} Final error (Euclidean norm) for the exact operator splitting solution~\eqref{phs.strang.1}-\eqref{phs.strang.3} (x), operator splitting with discrete gradient method~\eqref{dgm.1}-\eqref{dgm.2} ($\cdot$), multirate approach with nested integration (o) and multirate integration with higher order scheme ($\ast$).}
\end{figure}

Fig.~\ref{fig:solutions-2} shows the error plot for different operator splitting approaches applied to the two mass oscillator. All schemes show nicely an order two behaviour, they only differ in the error constant. Though both are of order two. the exact operator splitting 
is about $10^4$ times more accurate than operator splitting based on the discrete gradient method. The same applies to the multirate approaches (large step size with highly accurate scheme vs. nested integration with small time steps for the energy preserving part).

\subsection{A single mass-spring-damper chain}
\begin{figure}
    \centering
    \includegraphics[width=0.49\textwidth]{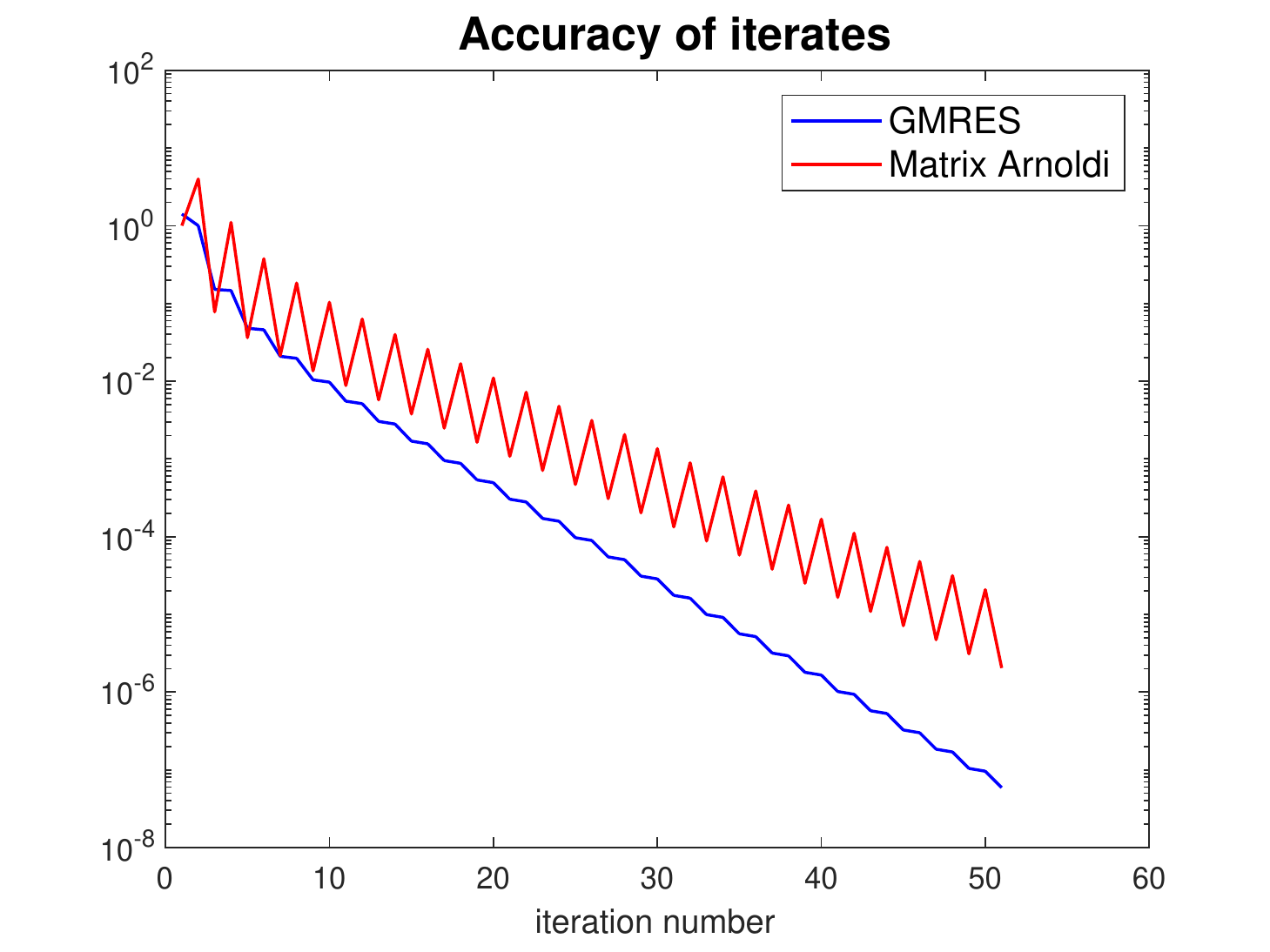}
    \hfill
    \includegraphics[width=0.49\textwidth]{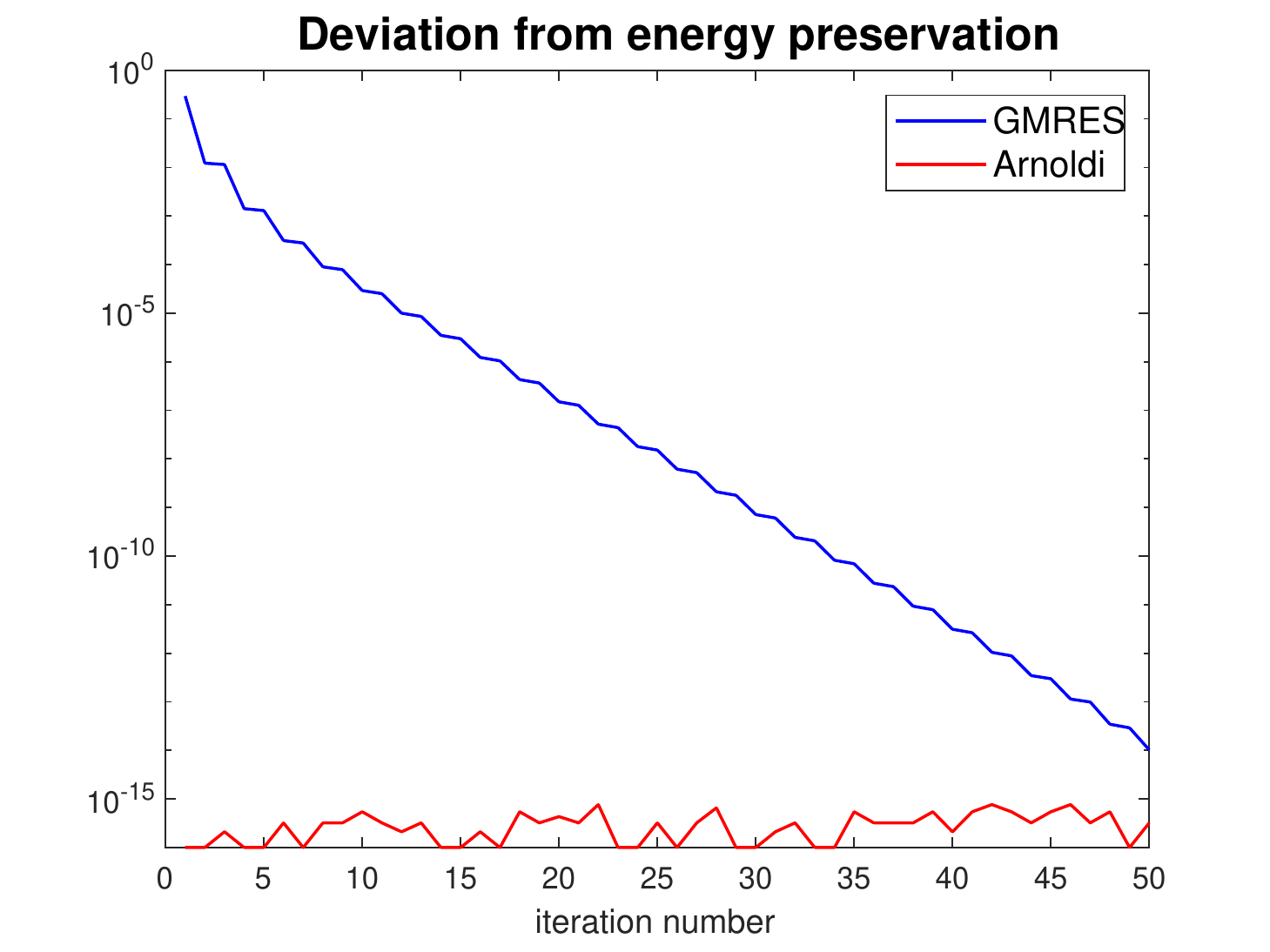}
    \caption{Comparison of GMRES and the matrix Arnoldi method. Left: 2-norm of the residual as a function of the iteration number. Right: $|1-\| x^{(k)}\|/\|\tilde{x}_0|\| \, |$ as a measure for the deviation in energy preservation.}
    \label{fig:compare_Arnoldi_GMRES}
\end{figure}

To highlight the structure preserving properties of the matrix Arnoldi approach for system~\eqref{eq:J_average_field} we use the single MSD chain from the Port Hamiltonian Benchmark System\footnote{\url{https://algopaul.github.io/PortHamiltonianBenchmarkSystems.jl}} with vanishing inputs where we chose the size parameters such that we obtain a dimension $N = 10,000$. For this example, the largest eigenvalue of $J$ is $10i$ so that we chose $h = 0.05 \cdot \frac{1}{10}$ as a step size which should sample the periods of all frequencies sufficiently well.
The results are given in Figure~\ref{fig:compare_Arnoldi_GMRES}

{The left part of Figure~\ref{fig:compare_Arnoldi_GMRES} reports the convergence history for two approaches to solve \eqref{eq:J_average_field}. The first approach just uses the GMRES method for the matrix $I-\frac{h}{2}\tilde{J}$, the second is the Arnoldi matrix function method. We report the size of the residuals 
$(I-\frac{h}{2}\tilde{J})x^{(k)}-(I+\frac{h}{2}\tilde{J})\tilde{x}_0$ for the iterate $x^{(k)}$ for both methods. We see that when measuring accuracy via the residual, the GMRES method is about 25\% faster than the matrix Arnoldi method. This corresponds to similar computational cost, since both methods require one matrix-vector multiplication with $\tilde{J}$ per iteration. The Arnoldi approach shows a decrease in quality in every other iteration. This can be attributed to the fact that for $k$ odd the projected matrix cannot well accomodate the symmetry of the spectrum of $\tilde{J}$ with respect to the real axis. 
}

{
The right part of Figure~\ref{fig:compare_Arnoldi_GMRES} shows that the matrix Arnoldi method does a perfect job in energy preservation, since the 2-norm of all its iterates differ from that of $\tilde{x}_0$ ,in a relative sense, by just machine precision. For the GMRES approximations, the violation of energy preservation is quite pronounced for the early iterates, and it becomes less as the iteration proceeds. 
}

\begin{figure}
    \centering
    \includegraphics[width=0.5\textwidth]{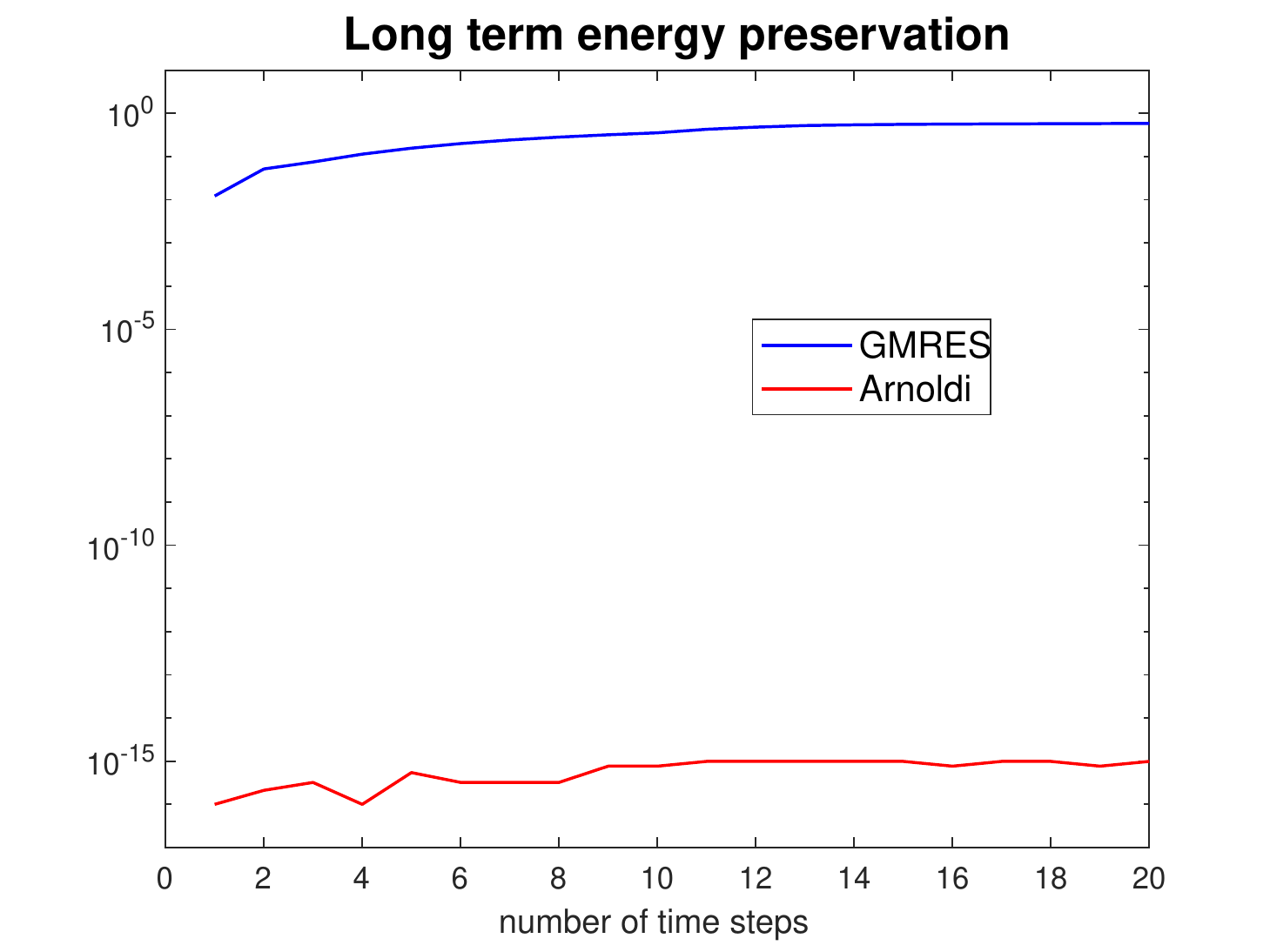}
    \caption{Deterioration of energy conservation after several time steps. We plot $|1-\| x_i \|/\|\tilde{x}_0\| \, |$, $x_i$ the numerical approximation for $x(t_0 + ih)$, as a measure for the deviation in energy preservation.}
    \label{fig:multirate}
\end{figure}

{
As a further illustration of the effects of the inherent structure preservation of the matrix Arnoldi approach, Figure~\ref{fig:multirate} reports a situation occurring within a multi-rate setting, see Section~\ref{sec:multirate}. We now do 20 consecutive steps of numerical integration, and in each step we stop the computation of the numerical approximation to the solution of \eqref{eq:J_average_field} once the residual is of size $h^2$. This choice is motivated by the fact that the integration scheme by itself has order 2. The figure again reports the quality of energy preservation as the relative difference of the 2-norms of the computed approximations for $x(t_0+ih), i=1,\ldots,20$ and the 2-norm of the initial value. We see that the GMRES approach now presents a very severe violation of energy preservation at later time points, whereas the matrix Arnoldi approach again preserves energy perfectly up to machine precision. We note that we required the same residual accuracy for both methods, which means that the matrix Arnoldi approach takes about 25\% more matrix vector products than GMRES.     
}

\section{Summary}
In this paper we have discussed how operator splitting methods at different levels can be used for the numerical simulation of port-Hamiltonian systems for both obtaining efficient and structure preserving methods: exact operator splitting based on Strang splitting, the use of discrete gradient schemes, exploiting the multirate behaviour in the splitting between structure-preserving and dissipative part, and structure preserving numerical solution of the respective linear systems by a matrix Arnoldi approach. 

Open questions for future research comprise, amongst others,  higher-order sche\-mes, tailored linear solvers for the dissipative part and generalization to port-Hamiltonian DAE systems.

\bibliography{ref}
\bibliographystyle{siam}

%
%
%



\end{document}